# Modeling branching effects on source-sink relationships of the cotton plant

Dong Li[1], Véronique Letort[2], Yan Guo[1], Philippe de Reffye[3, 4], Zhigang Zhan[1, *]

[1] Key Laboratory of Plant-Soil Interactions, Ministry of Education, College of Resources and Environment, China Agricultural University, Beijing 100193, China
[2] Ecole Centrale of Paris, Laboratoire de Mathématiques Appliquées aux Systèmes, F-92295 Châtenay-Malabry cedex, France
[3] Cirad-Amis, UMR AMAP, TA 40/01 Ave Agropolis, F-34398 Montpellier cedex 5, France
[4] INRIA-Futurs, EPI Digiplante, F-92295 Châtenay-Malabry cedex, France
[*] E-mail for correspondence: zhigang.zhan@cau.edu.cn

## Abstract

*Compared with classical process-based models, the functional-structural plant models provide more efficient tools to explore the impact of changes in plant structures on plant functioning. In this paper we investigated the effects of branches on the source-sink interaction for the cotton plant (Gossypium hirsutum L.) based on a two-treatment experiment conducted on cotton grown in the field: the single-stem plants and the plants with only two vegetative branches.*

*It was observed that the branched cotton had more organs for the whole plant but the organs on the trunk were smaller than those on the single-stem cotton. The phytomer production of the branches was four or five growth cycles delayed compared with the main stem. The organs on the trunk had similar dynamics of expansion for both treatments. Effects of branches were evaluated by using the functional-structural model GREENLAB. It allowed estimating the coefficients of sink strength to differentiate the biomass acquisition abilities of organs between different physiological ages. We found that the presence of the two vegetative branches increased the ground projection area of plant leaves and had led to slight changes on the directly measured parameters; the potential relative sink strengths of organs were found similar for the two treatments.*

## 1. Introduction

The functional-structural models simulate the architectural development and physiological functioning of plants. Most of them describe the source-sink relationships among the plant organs, see for instance L-PEACH for peach tree [1], LIGNUM for pine [2, 3].

GREENLAB is one of these models. It aims at simulating and predicting the growth and structural development of plants [4, 5]. The functional parameters that control the processes of biomass production and allocation are estimated in parallel from data of organ biomass. GREENLAB has been calibrated on some unbranched plants, such as maize [6-8], single-stem cotton [9], sunflower [10], and tomato [11]. The procedure of parameter estimation thus provides new insights about the plant functioning as it allows us to trace back the dynamics of the source-sink relationship within the plant [12]. A few studies have been undertaken to analyze the growth of branched plants with GREENLAB [13-15]. However, these studies were performed on trees grown in open field area that did not belong to any experimental setting. Therefore no pruning strategies were applied on these individuals.

However, it is of great interest to use functional-structural models to analyze the influence of plant structure on its functional processes. This is precisely the reason why such kind of models has been developed. In this paper we aim at understanding the impact of the structure of cotton plants (*Gossypium hirsutum* L.) on their physiology in terms of biomass production and allocation. Previous models developed to simulate the growth of cotton are mainly process-based models (e.g. GOSSYM [16] and OZCOT [17]), therefore the effect of the pruning strategy on the source-sink balance has not yet been studied. There are already some functional-structural models to simulate the branched plant, such as L-OZCOT for cotton [18], but little attention is put to the effect of branches.

In this paper, we present an experiment on pruned cotton with two different pruning treatments. In treatment 1, all branches were pruned (single-stem cotton); in treatment 2, only the two vegetative branches were kept. The goals of this study are (1) to observe the effects of vegetative branch presence on cotton development and growth; (2) to simulate the two treatments using the GREENLAB model; (3) to test whether the GREENLAB model can simulate the growth of cotton plants submitted to the two different treatments with a single set of parameters and, in the opposite case, to find which parameters are different.

## 2. Materials and methods

### 2.1. Field experiment

The field experiments were conducted at Quzhou experiment station (36°52' N, 115°1' E) in North China Plain, from 6 May to 25 August 2006. The cotton cultivar DP99B (*Gossypium hirsutum* L.) was used. The seeds were sown at a spacing of 0.8 m×0.8 m. Fertilizer inputs and irrigation were conducted so as to avoid any mineral or water limitations to plant growth.

Two pruned treatments were conducted in the course of plant growth. For the first treatment (T1), all the branches were removed immediately after their initiation [19]. The resulting plants were single-stemmed having only leaves kept on the main stem. For the other treatment (T2), all the branches were removed except the two vegetative branches on the 5th and 6th node were kept on the main stem, and both branches were pruned as the main stem. An example of the resulting plant was shown in Fig. 1.

Both non-destructive observations in the field and destructive measurement in the laboratory were conducted. We marked 10 single-stem cottons and every two days we measured the length of blade rib, petiole and internode to determine the expansion time of the organs.

Destructive samples were taken every two weeks. Each time, four plants per treatment were measured on dry weight and area for leaf blades, dry weight, length and diameter for internodes and petioles. Blade area was measured using a LI-COR 3100 leaf area meter (Lincoln, NB, USA). Roots were not taken into account. The thermal time was computed as the sum of mean daily air temperature minus a base temperature (15 ℃ in this study) after emergence. The climate data were from the climate station located near the experimental field.

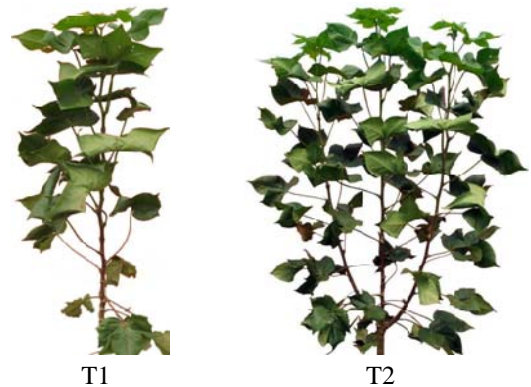

**Figure 1. Two treatments on cotton plants. T1 is the single-stem cotton having no branches. T2 is the cotton plant with two vegetative branches on the 5th and 6th node of the main stem.**

The measured organ dry weight and size data were input in a target file according to the topological structure which was subsequently used to estimate the model parameters. The target file format and optimization procedures for single-stem plants were described in [6]. For branched cotton, the position of the two vegetative branches was recorded and their biomass data were described the same way as for the main stem. The Digiplante software [20] was used to fit the GREENLAB model for the two treatments respectively. For each treatment, four-stage data were used as the target to fit the parameters.

### 2.2. The GREENLAB model and new features

The functional-structural GREENLAB model has been previously described in detail for non-branching plants [6-9]. There are some developments with the model for the branched plants [14, 15]. We will present here the main principles of this model for branched cotton.

We define Growth Cycle (GC), which is the time elapsed between emergence of two successive phytomers expressed in thermal time, as the simulation time step. Dates of events such as phytomer appearance, as well as biomass production and partitioning, leaf expansion time and longevity, are expressed in GC. At each GC except the first one, GREENLAB calculates the biomass production for the plant from the total functional blade area by an empirical nonlinear function representing intercepted light. It is assumed that the biomass for the first GC comes from the seed. Biomass partitioning among the growing organs is according to their relative sink strengths. Organs can expand during several GCs and their weight corresponds to the accumulative biomass they received. Their chronological age is defined as the time elapsed since their initiation, expressed in GC. Geometric dimensions of organs, e.g. blade area and internode length, are computed according to specific allometric rules that define organ shapes.

The pruned cotton produces simple plant architectures. Each phytomer consists of a few components: a node, an internode, and a leaf (except the first phytomer with two cotyledons). Nodes are difficult to split from internodes, so the studied organs are internodes, blades and petioles.

For cotton with two vegetative branches, the two branches are assumed to have the same physiological age (denoted by PA2) while the physiological age of the main stem is set to PA1 [21]. So there is only PA1 for single-stem cotton. In the model, organs of the same type (internodes, blades, petioles) have the same shape of sink variation along their life while they can have different coefficients of the potential sink strength according to their PA.

Biomass acquisition at GC($i$) is computed by Eqn. (1), based on the formalism of the Beer-Lambert's law for light interception in an homogeneous medium :

$$Q(i) = \frac{E(i) S_p}{R} \left(1 - \exp\left(-\frac{\sum_{PA=1}^{PA_m} \sum_{k=1}^{n_b(i,PA)} S_b(k,PA)}{S_p}\right)\right) \quad (1)$$

where $Q(i)$ is the biomass production during GC($i$); $E(i)$ is the average potential of biomass production during GC($i$) depending on environmental factors (e.g. light, temperature, wind) and it is assumed constant through all the growth cycles in this study; $PA$ is the physiological age of the axis and $PA_m$ is the largest $PA$ for the whole plant ($PA_m = 1$ for T1, and $PA_m = 2$ for T2); and $n_b$ is the number of functional leaves at the beginning of GC($i$) with certain $PA$ and it is relevant with blade functional time $T_f$; $S_b$ is the functional blade area (including leaves on the main stem and on the two branches); $R$ is an empirical resistance parameter of the plant; $S_p$ is the ground projection area of the leaves, which takes into account their inclination [9].

Biomass feeds all growing organs (sinks) through a globally shared reserve pool. The amount of biomass partitioned to individual organs is proportional to its relative sink strength which is a function of organ chronological age $k$ in terms of GCs, with respect to the type of organ $o$:

$$D_{o,PA}(k) = C_{o,PA} P_o f_o(k) \quad (2)$$

where the subscript $o$ denotes organ type: b, leaf blade; p, petiole; e, pith (In the model, internode is considered as an organ made up of pith, which elongates for a few GCs, and layers added during the secondary growth, which lasts during the whole growth period.); $D_{o,PA}$ is the sink strength of the organ; $P_o$ is the potential sink strength (for leaf blade $P_b=1$ is set as reference, therefore $P_o$ is a relative value); $C_{o,PA}$ is the coefficient of sink strength to differentiate the biomass acquisition abilities of organs between different physiological ages (for PA1, $C_{o,1}=1$ is set); $f_o(k)$ is an organ type specific function of sink variation in their developing period with one GC as time step, taking the beta-like form:

$$f_o(k) = \begin{cases} g_o(k)/M_o & (1 \leq k \leq T_o) \\ 0 & (k > T_o) \end{cases} \quad (3)$$

where $g_o(k) = \left(\frac{k-0.5}{T_o}\right)^{\alpha_o-1} \left(1 - \frac{k-0.5}{T_o}\right)^{\beta_o-1}$ and

$M_o = \max\{g_o(k) | 1 \leq k \leq T_o\}$; $\alpha_o$ and $\beta_o$ are parameters associated with the organ type $o$. The six parameters (three organ types, blade, pith and petiole) can be obtained simultaneously by fitting the model to measured data, as well as other parameters ($R$, $S_p$, $C_{o,PA}$ and $P_o$). $T_o$ is organ expansion time which can be obtained from the observed data.

The secondary growth for cotton covers all the existing internodes along the whole stem for each GC similarly to the process of ring growth of trees. In this study, virtual rings were defined to model secondary growth for cotton. The ring growth rate was found tightly dependent of the number of functional leaves left on the plants. Let $P_c$ be the potential sink strength of the ring compartment. Its sink strength, denoted by $D_c(i)$, is assumed to be proportional to the number of functional leaves, i.e. the sum of $n_b(i,PA)$ (the number of functional leaves during cycle GC($i$) of a certain PA; for branched cotton, the sum includes the functional blades on branches):

$$D_c(i) = P_c \sum_{PA=1}^{PA_m} n_b(i, PA) \quad (4)$$

One new ring is added to each phytomer at each growth cycle. The biomass allocated to internode ring is proportional to its length. The total biomass produced at GC(*i*) is distributed among all the organs according to the ratio of their relative sink to the total demand of the plant. And the organ biomass is cumulated GC by GC. At the end of GC(*i*), the weight of an organ with chronological age $k$ is calculated as:

$$Q_{o,PA}(i,k) = \frac{Q(i) D_{o,PA}(k)}{D_t(i)} + Q_{o,PA}(i-1,k-1) \quad (5)$$

where $D_t(i)$ is the total sink strength of all the organs including blades, petioles, piths and rings in GC(*i*).

As plant grows, the plant functional organs expand in mass according to their relative sink values and to the total biomass, which in turn depends on the functional blade area. This retroactive loop allows expressing the model in the form of dynamic system.

We set that internode pith and layers have the same density $\rho$ whatever their position in the plant is (branches or trunk). Volumes of internodes can be obtained from their weights.

The geometric shape for the pith is considered as a cylinder. For the volume of pith ($V_e$), the corresponding length $l$ and cross section $\sigma$ of the pith are then computed according to the relation with two allometric coefficients $b_e$ and $a_e$ [9, 19], *i.e.*,
$l = \sqrt{b_e} \ V_e^{(1+a_e)/2}$ and $\sigma = V_e / l$.

At each growth cycle, a new ring is added to the internode, so its diameter increases every GC.
To summarize, fifteen crop parameters that are not accessible from direct measurement are optimized by the fitting procedures. Four potential sink strengths and their coefficients for PA2 corresponding to the organ types (i.e., $P_b$ and $C_{b,2}$, $P_p$ and $C_{p,2}$, $P_e$ and $C_{e,2}$, $P_c$ for leaf blade, petiole, pith, and ring, respectively); For single-stem cotton there is only PA1, so the number of parameters is 12. Six coefficients of three beta function parameters are used to define the sink strength variation (i.e., $\alpha_b$ and $\beta_b$, $\alpha_p$ and $\beta_p$, $\alpha_e$ and $\beta_e$ for leaf blade, petiole and pith, respectively); two parameters (i.e., $R$ and $S_p$) for biomass acquisition. The allometric parameters including $\varepsilon$ (specific leaf weight), $b_e$ and $a_e$ (two shape coefficients of pith), as well as $Q_s$ (seed weight), $T_o$ (organ expansion time), $n_b$ (number of functional blades), $T_f$ (blade functional time) were directly measured and input in the model.

## 3. Results

### 3.1. Direct observation

The pruned cotton grew taller, had larger blades but less total biomass than the regular-growth cotton. And the blades of the pruned cotton were prone to twist when they grew larger (Fig. 1). The first and second branch was observed to emerge averagely six and five GCs after their bearing internode, respectively.

The organ size in two-vegetative branch cotton is smaller than single-stem cotton.

The directly measured parameters and their values were shown in Table 1. The organ expanding times were assumed the same for both treatments, i.e., 8 GCs for pith, and 15 GCs for blade and petiole. The functioning time was set to 22 for all the blades.

The blades in branched cotton were thinner than single-stem cotton but they were only about 5% difference between the two treatments for the main stem. The shapes of the pith for the two treatments were similar for the main stem and the branches were more different. And the internode densities of the two treatments were very similar, too.

**Table 1. The directly measured parameters and their values**

| Parameter | Treatment | | |
|---|---|---|---|
| | T1 | T2 | |
| | PA1 | PA1 | PA2 |
| $\varepsilon$ ($10^{-3}$ g cm$^{-2}$) | 6.33 | 5.98 | 5.59 |
| $b_e$ | 12.66 | 14.87 | 20.80 |
| $a_e$ | 0.047 | -0.028 | -0.0038 |
| $\rho$ (g cm$^{-3}$) | 0.26 | 0.25 | 0.25 |

$\varepsilon$ is specific leaf weight; $b_e$ and $a_e$ are coefficients of the allometric relationship for pith; $\rho$ is internode density. T1 denotes single-stem cotton and T2 denotes branched cotton.

### 3.2. Fitting results

The fitting target consisted of four-stage data of cotton for both treatments. For branched cotton, data were from the main stem and two branches. Fig. 2 showed fitting results on main stem for both treatments, and on the first branch for T2. Fitting on the second branch was similar as branch 1 (not shown). The model fitted the data well in both cases, especially for biomass of blade and petiole. The fitted internode biomass at the early stages was larger than the observation, and the fitted internode length was smaller.

The potential sink parameters and biomass production parameters optimized by the fitting process were shown in Table 2. The sink parameter values were nearly the same for the two treatments. The different parts were not shown (coefficients of PA2 for branched cotton, i.e., $C_{b,2}, C_{p,2}, C_{e,2}$ for blade, petiole and pith, respectively. Their values were less than 1.0). The value of $S_p$ for branched cotton was twice larger than that for single-stem cotton.

**Table 2. The optimized parameters and their values of the two treatments**

| Parameter | Treatment | |
|---|---|---|
| | T1 | T2 |
| $P_b$ | 1 | 1 |
| $P_p$ | 0.37 (0.05) | 0.36 (0.04) |
| $P_e$ | 0.36 (0.04) | 0.31 (0.04) |
| $P_c$ | 0.27 (0.02) | 0.23 (0.01) |
| $R$ | 192.9 (10.1) | 223.1 (6.8) |
| $S_p$ (cm$^2$) | 2209 (161) | 5475 (257) |

Parameter $P_b$, $P_p$, $P_e$, and $P_c$ are potential sink strength for blade, petiole, pith, and ring, respectively ($P_b$ is set to 1 as reference); $R$ is coefficient for plant resistance to transpiration; $S_p$ is light interception parameter for Beer–Lambert's law. Data in brackets were standard errors of the estimated parameter values; T1 and T2 are specified as in Table 1.

Compared with single-stem cotton, all the sink variation functions for branched cotton were asymmetrically shrunk, and the time when the maximum growth rate occurred for both blade and petiole were delayed (Fig.3). The branches influence on sink variation for pith was the least.

The 3D representation of cotton plants at four stages was also generated (Fig.4).

## 4. Discussion and conclusions

The presence of branches not only decreased organ (blade, petiole and internode) size on the trunk (Fig. 2), but also changed the organ expansion function (Fig. 3). Branches can change plant trunk in a few aspects, for example, the specific leaf weight of main stem for cotton of two vegetative branches was smaller than single-stem cotton; the internode allometric parameters were different. The size and shape of organs on branches were different from those on the trunk, and accordingly parameters varied in GREENLAB with PAs. The model can simulate the growth of both treatments well, so the fitted values of parameters may reflect the influence of branch presence.

The potential sink strength ($P_o$) was very close for the two treatments (Table 2). Because the beta function is very flexible, the sink variation is similar between the two treatments (Fig.3) despite of the differences of the parameters, especially when it shifts one cycle since the value of the initial stage is very small. In fact, we can see that the fitting results of two-vegetative branch cotton can be improved and maybe more similar organ expansion function can be obtained as the single-stem cotton. It indicates that there is no significant change in their expansion type, just as shown in the results part. So for the sink parameters, we can use the same values to simulate the trunk growth for the two treatments.

The ground projection area ($S_p$) was acquired by the fitting optimization procedure. The value range of this parameter can also be estimated from the observed data [12]. For branched cotton, the projection area was made up of contributions of leaves both on trunk and on branches. So the existence of branches may increase the value of $S_p$. Because it is difficult to separate the pith from the ring, the allometry of pith was not precisely obtained. Based on these parameters, the fitted internode length was a bit smaller than measured. The allometric parameters may be improved by optimization procedures.

For a first approximation, we used a constant projection area ($S_p$) for all plant growth stages. A further study is undergoing to take into account the dynamics of this parameter.

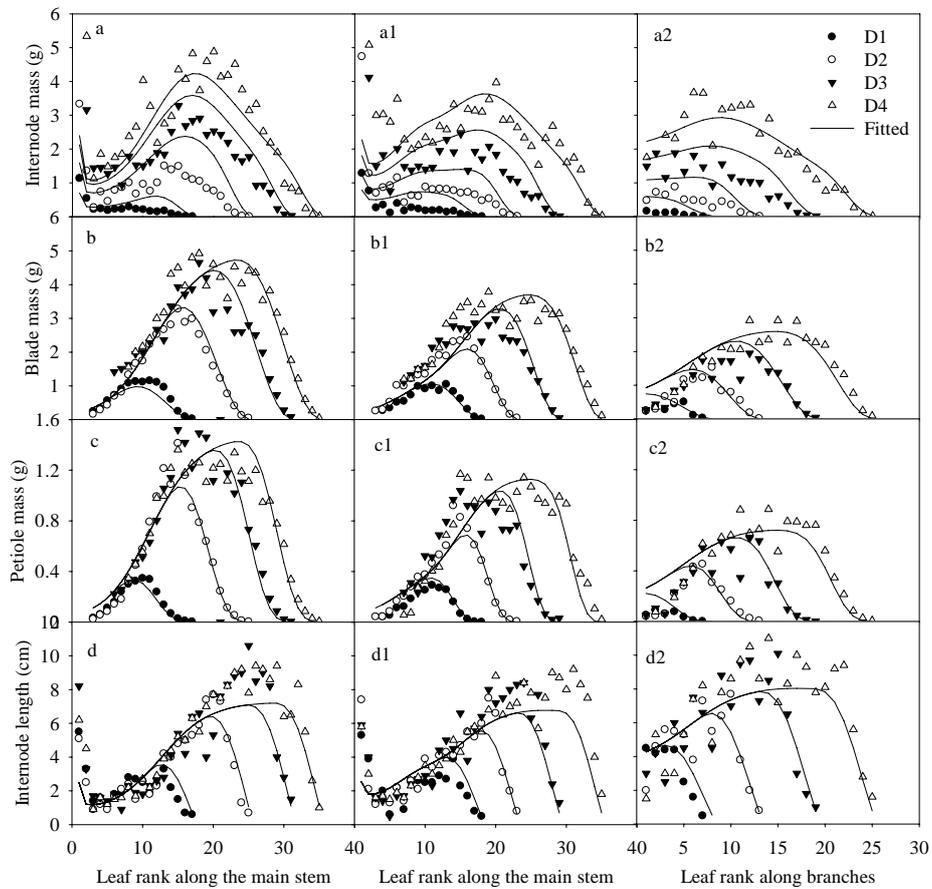

**Figure 2. Fitting GREENLAB model to four-stage data of cotton plants. Letters a, b, c and d denote single-stem cotton; a1, b1, c1, d1 denote main stem of branched cotton; a2, b2, c2 and d2 denote the first branch. The target items include biomass and length of internode (a, d), biomass of leaf blade (b) and petiole (c). D1, D2, D3, and D4 denote four growth stages. The observed data are the samples which are nearest to the mean of four replications.**

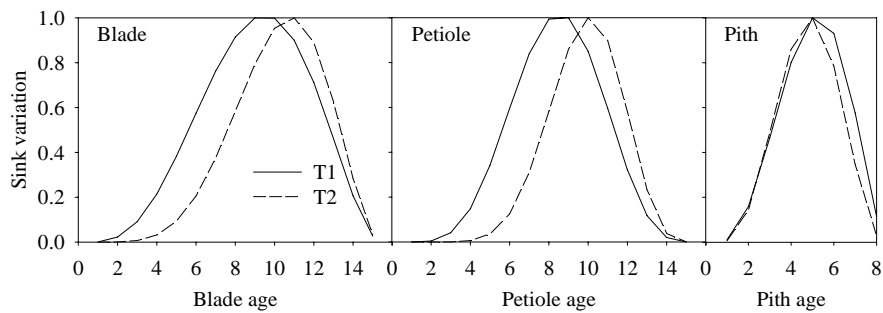

**Figure 3. Comparison of organ sink variation functions for the two treatments. T1 and T2 are specified as in Table 1.**

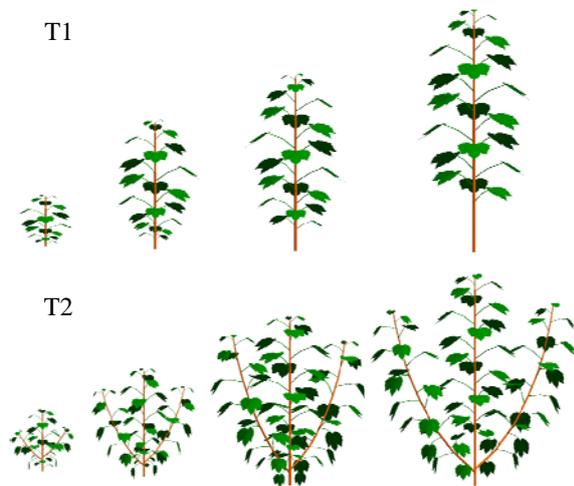

**Figure 4.** The 3D representation of cotton plants at four stages of 18, 25, 30, 35 GCs, respectively. T1 and T2 are specified as in Table 1.

## Acknowledgements

This research has been partly supported by the Hi-Tech Research and Development (863) Program of China (2006AA10Z229). The authors are grateful to Chuntang Gong, Yufei Li and Zhifei Li for their kind help on the experiment measurements, and to Paul-Henry Cournède for allowing the use of the Digiplante software developped at the Laboratory of Applied Mathematics in Ecole Centrale of Paris.